\documentclass[reprint,aps,prx,onecolumn,superscriptaddress]{revtex4-2}
\usepackage{hyperref}
\usepackage{amsmath,amsthm,amssymb,amsthm}
\usepackage{braket,physics}
\usepackage{enumitem}
\usepackage{color}
\usepackage{cancel}
\usepackage{tikz}
\usepackage{graphicx}
\usepackage{subcaption}
\usepackage{booktabs,tabularx,multirow}

\begin{document}

\title{Classical Combinatorial Optimization Scaling for Random Ising Models on 2D Heavy-Hex Graphs}

\author{Elijah Pelofske}
\email[]{epelofske@lanl.gov}
\affiliation{Information Systems \& Modeling, Los Alamos National Laboratory, USA}

\author{Andreas Bärtschi}
\affiliation{Information Sciences, Los Alamos National Laboratory, USA}

\author{Stephan Eidenbenz}
\affiliation{Information Sciences, Los Alamos National Laboratory, USA}

\begin{abstract}
Motivated by near term quantum computing hardware limitations, combinatorial optimization problems that can be addressed by current quantum algorithms and noisy hardware with little or no overhead are used to probe capabilities of quantum algorithms such as the Quantum Approximate Optimization Algorithm (QAOA). 
In this study, a specific class of near term quantum computing hardware defined combinatorial optimization problems, Ising models on heavy-hex graphs both with and without geometrically local cubic terms, are examined for their classical computational hardness via empirical computation time scaling quantification. Specifically the Time-to-Solution metric using the classical heuristic simulated annealing is measured for finding optimal variable assignments (ground states), as well as the time required for the optimization software Gurobi to find an optimal variable assignment. 
Because of the sparsity of these Ising models, the classical algorithms are able to find optimal solutions efficiently even for large instances (i.e., $100,000$ spin variables). The Ising models both with and without geometrically local cubic terms exhibit average-case linear-time or weakly quadratic scaling when solved exactly using Gurobi, and the Ising models with no cubic terms show evidence of exponential-time Time-to-Solution scaling when sampled using simulated annealing. These findings point to the necessity of developing and testing more complex, namely more densely connected, optimization problems in order for quantum computing to ever have a practical advantage over classical computing. Our results are another illustration that different classical algorithms can indeed have exponentially different running times, thus making the identification of the best practical classical technique important in any quantum computing vs. classical computing comparison.

\end{abstract}

\maketitle

\vspace{-0.5cm}

\section{Introduction}
\label{section:Introduction}

The Quantum Approximate Optimization Algorithm (QAOA) is a quantum heuristic for solving combinatorial optimization problems \cite{farhi2014quantumapproximateoptimizationalgorithm, Hadfield_2019}. QAOA is a promising quantum algorithm as it can be executed on noisy current quantum computers that are not error-corrected and some signal can still be measured despite the noise. There is evidence that QAOA could be competitive with classical algorithms \cite{Harrigan_2021, Shaydulin_2024, PhysRevA.104.052419, farhi2020quantumapproximateoptimizationalgorithm1, farhi2020quantumapproximateoptimizationalgorithm2, Farhi2022quantumapproximate, boulebnane2022solvingbooleansatisfiabilityproblems, montanaro2024quantumspeedupssolvingnearsymmetric}. 
QAOA and similar quantum algorithms that solve combinatorial optimization problems are of considerable interest because combinatorial optimization problems are a fundamentally challenging and important computational task in many fields including information processing -- and therefore if quantum mechanical effects can be utilized to efficiently solve such problems, then this would be an excellent use for quantum computing. 
However, studying quantum algorithms performance such as QAOA is inherently challenging because simulations of the algorithm requires either classical simulation methods, which are severely limited either due to the exponential classical computational complexity of simulating many-body quantum systems, or requires using current state of the art quantum computers which are heavily affected by limited qubit coherence times and various sources of control error. Thus, generally simulations to study QAOA performance are limited to either i) small problem sizes, or ii) intermediate-scale hardware with noise, or both. This limitation of studying quantum algorithms has motivated several different studies, including many different computational tasks not necessarily only for combinatorial optimization, which define problem types that can be implemented on NISQ (Noisy Intermediate-Scale Quantum) \cite{Preskill_2018} hardware with little or no overhead with respect to compute time or qubit use \cite{tasseff2022emergingpotentialquantumannealing, Andrist_2023, QAOA_QA_127, pelofske2023short, pelofske2023scalingwholechipqaoahigherorder, PhysRevX.4.021008, kim2023evidence, kechedzhi2023effective, begušić2023fast, tindall2023efficient, begušić2023fast_2, liao2023simulation, rudolph2023classical, patra2023efficient}. These types of studies allow various quantum algorithms to be investigated in a regime where some amount of signal can be obtained from the noisy quantum hardware, but these problem instances are typically not very computationally hard (with a few exceptions, such as ref.~\cite{tasseff2022emergingpotentialquantumannealing}) or particularly computationally relevant to a specific industry combinatorial optimization problem.

In this study, one of these types of NISQ-friendly problem instances is examined for its computational hardness using classical algorithms. We focus on two approaches: the commercial optimization software Gurobi \cite{gurobi} and the heuristic Simulated Annealing \cite{kirkpatrick1983optimization}. Simulated annealing is very general-purpose, meaning it can work quite well on a variety of different types of optimization problems, and therefore is typically considered as a performance baseline. 
Gurobi on the other hand serves as a powerful mathematical optimization solver tool, which in general we expect is closer to being state of the art for general purpose optimization.
We note that simulated annealing is not a state of the art heuristic solver of optimization problems -- there exist algorithms which perform better than simulated annealing (especially when considering specific types of optimization problems), such as Simulated Quantum Annealing~\cite{simulated_quantum_annealing}, Tabu Search~\cite{tabu_search}, Parallel Tempering~\cite{parallel_tempering_ICM} and others, which for example are benchmarked in ref.~\cite{tasseff2022emergingpotentialquantumannealing}.

The problem instances that we consider in this study are a type of discrete optimization problem, specifically random Ising models. The computational goal is to find the global minimum energy of the Ising model, or equivalently an optimal decision variable assignment. 
Specifically, this study examines the classical computational hardness of ground-state finding for two classes of random spin glass Ising models defined on heavy-hex graphs, those with only linear and quadratic polynomial terms and those with additional cubic (three body) terms, were introduced and defined in refs.~\cite{QAOA_QA_127, pelofske2023short}. This class of Ising models is designed to be extremely hardware compatible for current IBM Quantum superconducting qubit processors which have two-qubit gate operations that define a sparse graph structure that is generally known as heavy-hex \cite{Chamberland_2020, PRXQuantum.5.040334, benito2025comparativestudyquantumerror}. These heavy-hex Ising models, and in particular the implementation of the short-depth QAOA circuits that can sample these Ising models, have been subsequently been used in several studies \cite{barron2023provableboundsnoisefreeexpectation, pelofske2023scalingwholechipqaoahigherorder, sachdeva2024quantumoptimizationusing127qubit, mcgeoch2024commentcomparingoptimizationdwave}.

\section{Methods}
\label{section:methods}

\begin{figure}[h]
    \centering
    \includegraphics[width=0.32\linewidth]{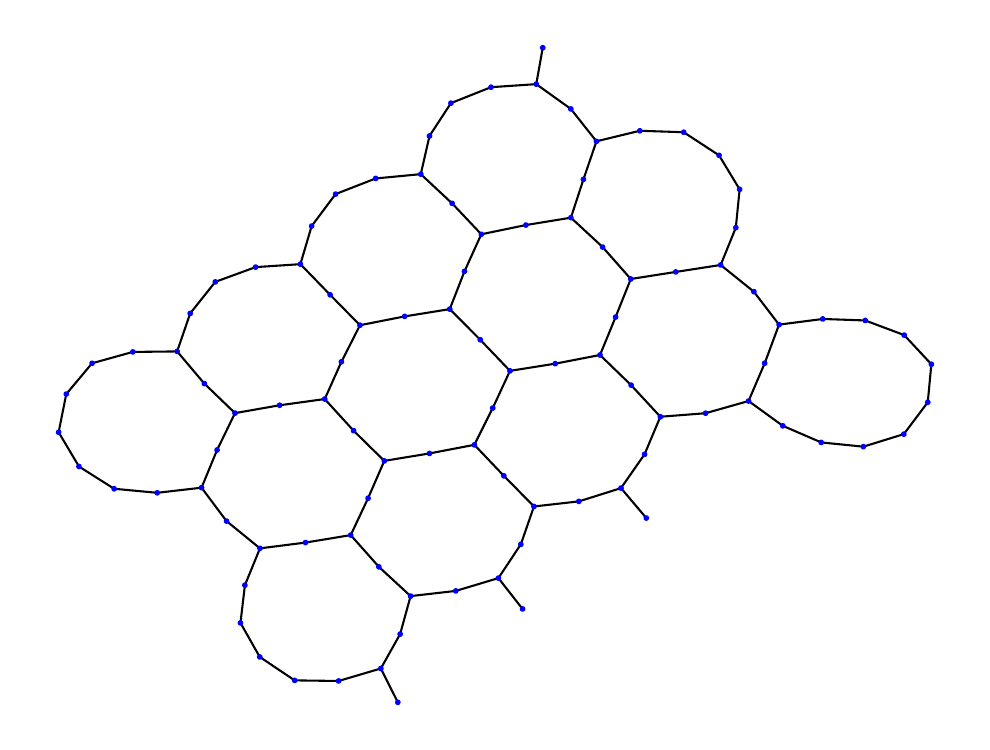}
    \includegraphics[width=0.32\linewidth]{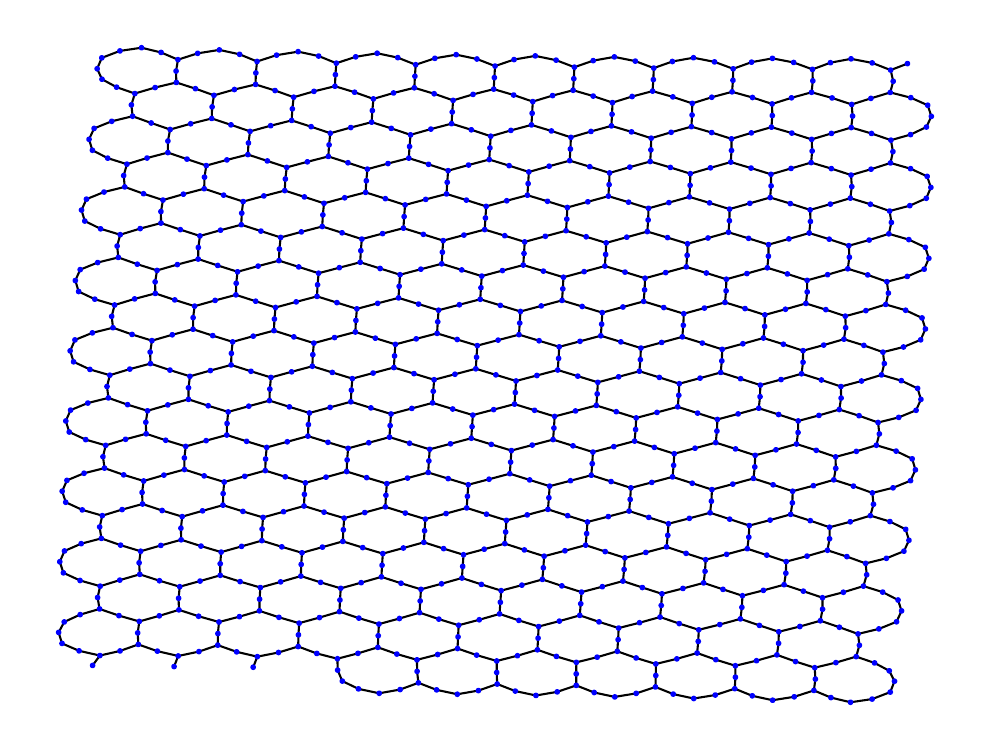}
    \includegraphics[width=0.32\linewidth]{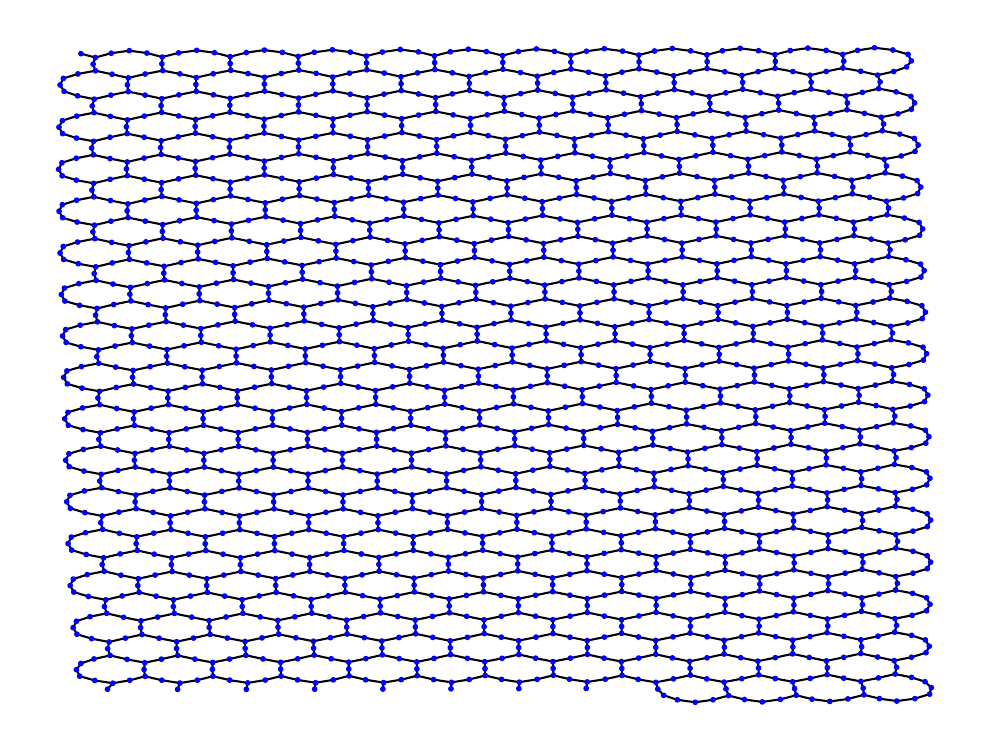}
    \caption{Examples of heavy-hex graph structures with 100 nodes (left), 1000 nodes (middle), and 2000 nodes (right). All three of these were created with the custom graph generator we created to produce approximately square lattices that follow the IBM heavy-hex hardware structure (although note that none of these structures are exactly the same as any current IBM quantum processor). }
    \label{fig:example_heavy_hex_graphs}
\end{figure}

First, we created a custom generator that produces heavy-hex graphs -- the structure of these graphs is a tileable lattice \cite{Chamberland_2020}. This allows us to create graphs with an arbitrary number of nodes. We focus specifically on 2-dimensional heavy-hex structures that are approximately square (approximately the same number of hexagonal unit cells on all sides of the 2D lattice) because this follows the existing structures of current IBM quantum processors \cite{Chamberland_2020}. Figure \ref{fig:example_heavy_hex_graphs} shows three of these generated heavy-hex graphs. Based on these graph structures the Ising models defined on those graphs can be created. The heavy-hex Ising model generator uses the Networkx package \cite{hagberg2008exploring}. 

The heavy-hex graph spin glass Ising models were generated using the same methods outlined in refs.~\cite{QAOA_QA_127, pelofske2023short, pelofske2023scalingwholechipqaoahigherorder}, but applied to significantly larger problem sizes. In the context of the field of combinatorial optimization, these Ising models are can be classified as \emph{minimization} and \emph{unconstrained} combinatorial optimization problems (as opposed to, for example, a type of constrained optimization problem). 

All of the decision variables in the Ising models are spins $z_i \in \{+1, -1\}$. Each heavy-hex graph can be described as a set of edges $E$ and vertices $V$; $G=(V,E)$. $|V|$ is the number of nodes which we will typically denote as $n$ to define the system size. Given a vector of spins $z = (z_1, \ldots, z_{n}) \in \{ +1,-1\}^n$, the cost function is defined as

\begin{equation}
    C(z) = \sum_{v \in V} d_v \cdot z_v + \sum_{(i,j) \in E} d_{i,j} \cdot  z_i \cdot z_j + \sum_{l \in W} d_{l,n_1(l),n_2(l)} \cdot z_l \cdot z_{n_1(l)} \cdot z_{n_2(l)},
    \label{equation:problem_instance}
\end{equation}

Typically we will refer to the cost value for a given sample as the energy, where the goal is to minimize the energy. Every heavy-hex graph is a bipartite graph that has vertices $V = \{1,\ldots,n\}$ and can be uniquely bipartitioned as $V = V_2 \sqcup V_3$ with $E\subset V_2\times V_3$, where $V_i$ has vertices of greatest degree $i$. $W$ is the set of vertices $l$ in $V_2$ that have degree $2$ and the neighbors are $n_1(l)$ and $n_2(l)$. $d_v$, $d_{i,j}$, and $d_{l,n_1(l),n_2(l)}$ are the linear, quadratic and cubic coefficients. The Ising model coefficients $d$ are randomly drawn from $\{+1, -1\}$ with probability $0.5$. This type of Ising model could be considered a type of sparse $\pm J$ model~\cite{PhysRevLett.92.117202, PhysRevLett.107.047203}. In this study we examine the computational hardness of ground-state finding for the Ising models with and without the geometrically local cubic terms.

The primary classical solver that will be used in this study is the commercial Gurobi optimization software \cite{gurobi} - which is not a single algorithm, but rather a piece of software that uses a variety of techniques \cite{achterberg2020presolve, achterberg1u2014multi, bixby1999mip, gu2000sequence, branch_cut_numerics}, such as branch-and-bound, and is designed to deterministically find an optimal solution, with a provable optimality gap so that it can be guaranteed there exists no better variable assignment (given sufficient compute time). Gurobi is not restricted to solving only discrete optimization problems, but here we use it to solve these discrete optimization problems where the variable states or either $-1$ or $+1$. 

Gurobi does not directly handle higher order terms in the objective function polynomial where the variable states are strictly spins ${+1, -1}$, and therefore they must be structured so as to be compatible with the software. We utilize two different formulation techniques to address the cubic terms and obtain global optimal solutions of the optimization problems:

\begin{enumerate}[noitemsep]
    \item The cubic $ZZZ$ terms in the higher order Ising models are order-reduced by introducing auxiliary variables, and quadratic interactions with those auxiliary variables, that preserve the optimal solution(s) to the original problem. This order reduction is performed using the heuristic order reduction method in the Python 3 module \texttt{dimod}\footnote{\url{https://github.com/dwavesystems/dimod}}. This order reduction is deterministic meaning that the resulting linear and quadratic terms in the polynomial are always the same given the same set of parameters (and there is no random seed). There is a penalty weight parameter that must be tuned such that the optimal (globally minimum) solutions are preserved in the order-reduced problem. Setting this penalty weight optimally is computationally hard in general. In this case, we set the weights to a large positive value equal to the number of nodes in the heavy-hex graph since we have found that in general as long as the weight is sufficiently greater than the largest magnitude coefficient in the original polynomial, the optimal solutions are preserved. There are very likely alternate order reduction techniques which do not require using a large penalty weight, however, order reduction methods are not the specific focus of this study and in particular this method empirically guarantees optimality in the tests we performed. This order reduced polynomial contains only linear and quadratic terms - and in Gurobi is formulated as a binary Integer Quadratic Program (IQP), where the binary variables (which are the variable types that Gurobi handles) are converted into spins via $2x_i-1$. Importantly, this technique introduces a large overhead of auxiliary variables, and in particular is not heavily optimized. More complex and tailored order reduction could certainly be applied to these types of higher order terms, but here we utilize general purpose heuristics that work efficiently. 
    \item Formulate the polynomial as a binary Integer Linear Program (ILP), where the variable states of the optimization problem are first converted into binary states. In order to convert variable states from spins into binary, with the higher order terms, the Python 3 module \texttt{dimod} is also used, which introduces additional higher order terms. The higher order terms are then \emph{linearized} \cite{ForresterHunt-Isaak2020, Padberg1989, Glover1975, SheraliSmith2007} into valid inequality constraint representations that do not use the higher order variable interactions directly, by using lower dimensional variable inequalities. An example is given below. All higher order terms can then be represented as a binary Integer Linear Program. This ILP formulation follows the same basic methods used in refs. \cite{tasseff2022emergingpotentialquantumannealing, dash2013notequboinstancesdefined}. 
\end{enumerate}

Here we show a small example of linearization of a cubic term, which is used in the binary ILP formulation of the problem instances. If we want to linearize the following cubic term,

$$K = x_i x_j x_k,$$ where $x_i, x_j, x_k$ are binary variables. Then we can specify a logically equivalent formulation using the following linear constraints:

\begin{align*}
K & \leq x_i \\
K & \leq x_j \\
K & \leq x_k \\
K & \geq x_i + x_j + x_k - 2
\end{align*}

The energy of the solution vectors produced by Gurobi when these two different methods are used, are compared against each other (for a given problem instance) to ensure that the two methods agree. In practice we found that these methods always agreed - however it is not clear whether the penalty weight used for the first method will always ensure an optimal solution is correctly encoded in the order-reduced polynomial for arbitrary polynomials.

In addition to being solved exactly using Gurobi, these Ising models are also sampled using an existing C++ implementation of simulated annealing that is available as a Python 3 library \texttt{dwave-neal}\footnote{\url{https://github.com/dwavesystems/dwave-neal}}. The simulated annealing algorithm performs a fixed set of Metropolis-Hastings spin updates, which we vary to be either $100$ or $1000$ where each update uses a schedule of increasing $\beta$ values. The default \emph{geometric} simulated annealing schedule is used, which defines the values of $\beta$ at each step of the schedule. $\beta$ is the inverse temperature of the Boltzmann distribution $\frac{1}{k_BT}$, $T$ is the thermodynamic temperature, and $k_B$ is the Boltzmann constant \cite{kirkpatrick1983optimization}. The total number of simulated annealing simulations performed for each Ising model, where each simulation generates a single low-energy sample of the optimization problem, is fixed to be $50,000$. The goal was to use a fixed number of samples that was reasonably large so that finite sampling does not dominate the reported statistics, but also computationally feasible to carry out the simulations in a reasonable amount of time. We use simulated annealing sampling strictly on the Ising models with no cubic terms, so as to avoid the increased computational overhead that comes from order reduction of the cubic terms (this particular simulated annealing solver works only with linear and quadratic terms). The computation quality of the simulated annealing sampling is measured using the well established heuristic probabilistic sampler based compute time metric \textit{Time-to-solution} (TTS) \cite{R_nnow_2014}:

\begin{equation}
    TTS = \frac{T_{\text{CPU}}}{N_s}  \frac{\log (1-0.99)}{\log (1-p)},
\end{equation}

where $p$ denotes the measured optimal solution sampling success rate (sampling rate of the optimal cost value, not any particular variable configuration). The optimal solution cost value is found using Gurobi. $N_s$ is the total number of samples obtained, $T_{\text{CPU}}$ is the classical CPU time required to generate the full set of $N_s$ samples measured as Python single-process time. TTS was defined in the context of quantifying heuristic probabilistic sampling quantum and classical algorithm runtimes, and serves as a good runtime complexity measure in the absence of guarantees of optimality \cite{R_nnow_2014, PhysRevLett.115.230501, Pearson_2020}. In this case TTS is defined as the expected compute time required to sample an optimal solution of the optimization problem at least once, with $99\%$ confidence. The optimal solution sampling success rate $p$ of finding the ground-state energy in this case is measured using the exact solutions obtained from Gurobi for the same problem instances. 

\begin{figure}[ht!]
    \centering
    \includegraphics[width=0.45\linewidth]{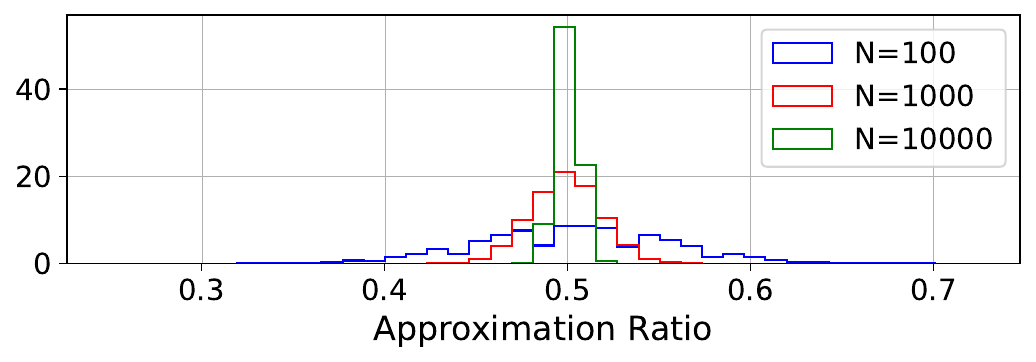}
    \includegraphics[width=0.45\linewidth]{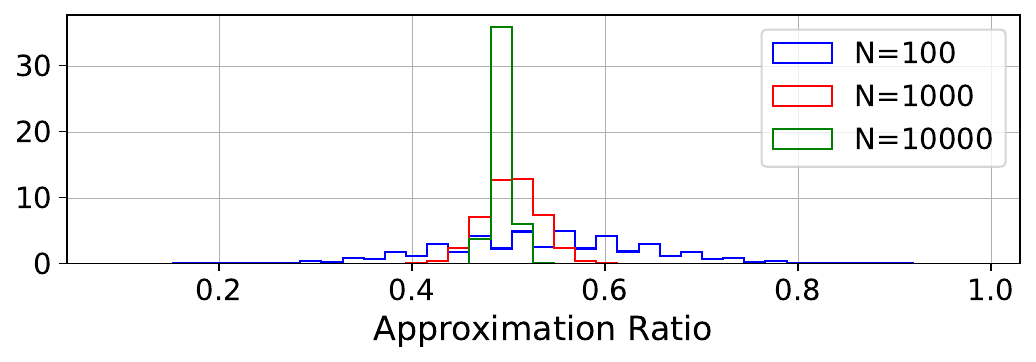}
    \caption{Normalized random sampling energy distributions, shown as density histograms (y-axis) for increasing problem size instances with cubic terms (left) and no cubic terms (right). The random samples are all generated with $p=0.5$ of both $-1$ and $-1$ states. Each distribution is created from $2\cdot10^{5}$ samples. Each distribution shows samples from a single problem instance. }
    \label{fig:random_sampling_energies}
\end{figure}

Gurobi Optimizer \cite{gurobi} version 11.0.3 is used for all Gurobi computations. The compute hardware used for simulated annealing simulations as well as Gurobi computations is a Red Hat Linux node with Intel(R) Xeon(R) CPU E5-2695 v4 2.10GHz. All Gurobi simulations use a time limit of $100,000$ seconds, a single thread, and a MIP gap of $1\mathrm{e}{-8}$. Gurobi was run with the intention of converging on an optimal solution, and each Gurobi run was checked for this condition where the optimality gap had closed and therefore the solution that we have is a verified optimal solution. At some point, for sufficiently large problems, the Gurobi simulation time would exceed this time limit that was specified, but for all simulations reported on in this study, each individual run completed with a converged optimality gap within the specific time limit. All binary solutions found by Gurobi are adapted to spins using the mapping of $1 \mapsto 1$, $0 \mapsto -1$. Gurobi can be used to find multiple feasible optimal solutions, if they exist, however in this case we use Gurobi to only find a single optimal solution for each problem instance. We use only a single thread in order to make the compute time scaling measurements very clear, and in particular to avoid multi-threaded CPU time measurements.

All curve fits are performed using non linear least-squares fitting using the scipy \cite{2020SciPy-NMeth} Python 3 library. Bounds are used to restrict the least squares fitting to not set the fitting variables too close to $0$ with a minimum coefficient precision limit of $10^{-7}$, for example in the exponential and quadratic cases, which would effectively be replicating a linear function scaling for small system sizes. All fitting coefficients are also forced to be non-negative. The linear function is fit to $y=a x + b$, log function is fit to $y=a \ln{x} + b$, quadratic function is fit to $y=a x^2 + b \cdot x + c$, and exponential function is fit to $y=e^{x a} \cdot b$.

\begin{figure}[ht!]
    \centering
    \includegraphics[width=0.58\linewidth]{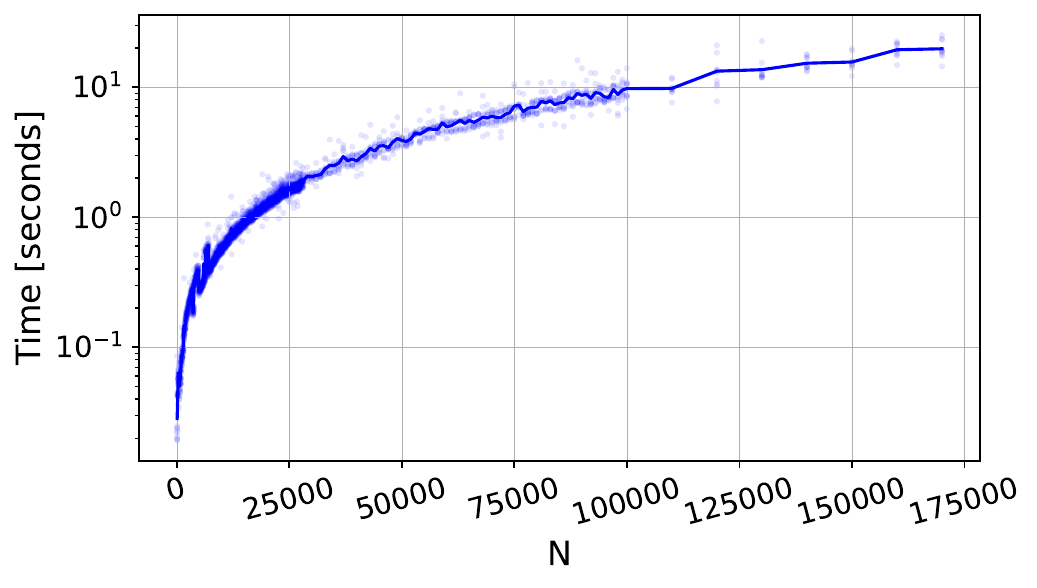}
    \caption{Gurobi runtime scaling (log scale y-axis) as a function of Ising model problem size (x-axis), where the instance contain no higher order terms. The x-axis (problem size) is exactly the number of nodes in the heavy-hex graph. For each system size, $10$ random optimization problems are solved and the runtimes are plotted as individual points; the mean time of each set of $10$ instances is plotted as a connected line. The problem formulation used in these simulations was binary IQP. }
    \label{fig:quadratic_Gurobi_time}
\end{figure}

\begin{figure}[ht!]
    \centering
    \includegraphics[width=0.49\linewidth]{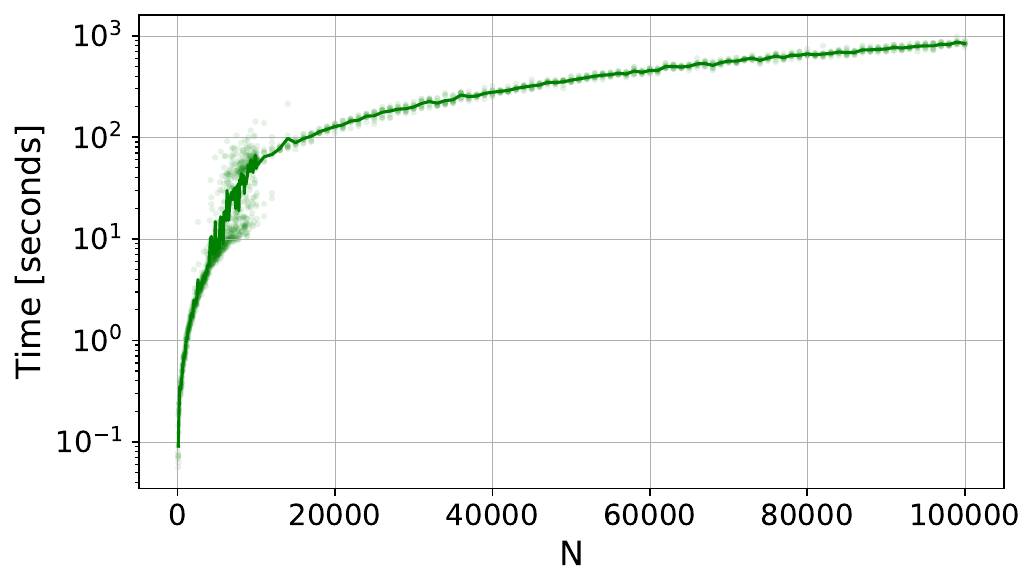}
    \includegraphics[width=0.49\linewidth]{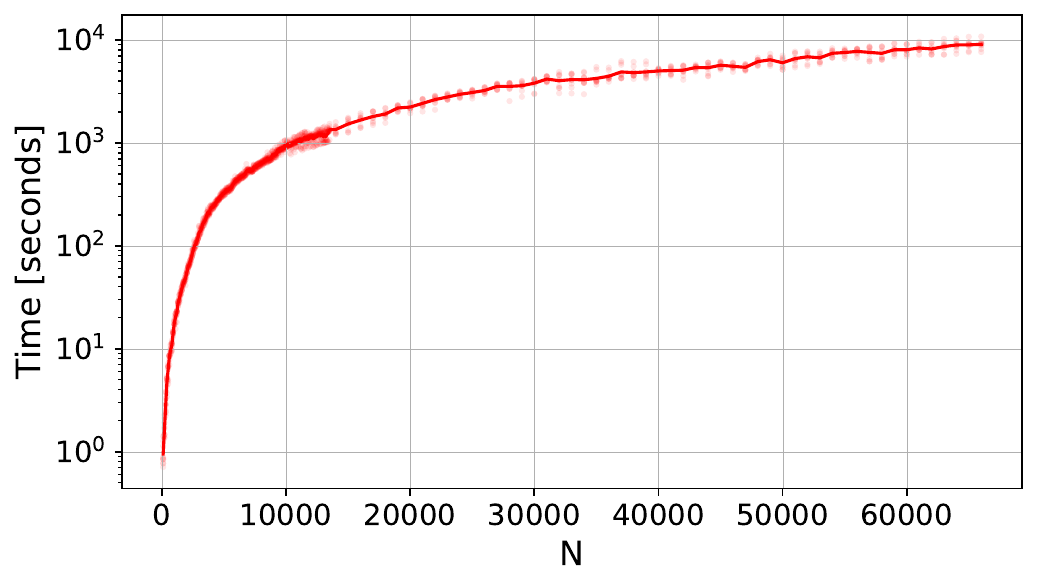}
    \caption{Gurobi runtime scaling (log scale y-axis) as a function of Ising model problem size (x-axis) for the instances with geometrically local cubic terms. The x-axis (problem size) is exactly the number of nodes in the heavy-hex graph. Left plot shows the runtime scaling for the binary ILP method, and the right plot shows the runtime scaling for the binary IQP order reduction method. These two methods of handling the cubic terms are applied to exactly the same problem instance. At each system size N, $10$ instances are solved by each method, and the individual runtimes are shown as individual points. Notice that the binary IQP method (right) requires approximately an order of magnitude more compute time than the binary ILP method (left) at $60,000$ variable problem sizes. }
    \label{fig:cubic_Gurobi_time}
\end{figure}

\section{Results}
\label{section:results}

\begin{table}[ht!]
{       \addtolength{\tabcolsep}{2pt}
        \newcolumntype{R}{>{\raggedleft\arraybackslash}X}  
        \begin{tabularx}{\textwidth}{@{}	l	l	l	R@{\hspace{4em}}	l	r	@{}}
            \toprule
            Problem type	&	Simulation	&	Function type	&	Function fit	&	Precise coefficient fits	&	RMSE	\\
            \midrule
            Quadratic	&	SA TTS, $1000$ updates
            			&	log			&	$37.2\,\ln x$						&	$a=37.2,\; b=10^{-6}$						&	608.9	\\
            			&				&	linear		&	$0.123\,x$						&	$a=0.123,\; b=10^{-7}$						&	534.3	\\
            			&				&	quadratic	&	$3.2\cdot10^{-5}\,x^{2}$				&	$a=3.2\cdot10^{-5},\; b=c=10^{-7}$			&	464.2	\\
            			&				&	\textcolor{red}{\bf{exponential}}
            										&	\textcolor{red}{$0.022\,e^{0.002\,x}$}
            										&	$a=0.002,\; b=0.022$							&	\textcolor{red}{235.7}	\\
            \cmidrule{2-6}
            			&	SA TTS, $100$ updates
            			&	log			&	$4.01\,\ln x$						&	$a=4.01,\; b=10^{-6}$						&	55.3	\\
            			&				&	linear		&	$0.0385\,x$						&	$a=0.0385,\; b=10^{-7}$						&	48.7	\\
            			&				&	quadratic	&	$3.6\cdot10^{-5}\,x^{2}$				&	$a=3.6\cdot10^{-5},\; b=c=10^{-7}$			&	41.8	\\
            			&				&	\textcolor{red}{\bf{exponential}}
            										&	\textcolor{red}{$0.000128\; e^{0.00898\,x}$}
            										&	$a=0.00898,\; b=0.000128$						&	\textcolor{red}{4.5}	\\
            \cmidrule{2-6}
            			&	Gurobi IQP
            			&	log			&	$0.24\,\ln x$						&	$a=0.24,\; b=10^{-6}$						&	2.7	\\
            			&				&	\textcolor{red}{\bf{linear}}
            										&	\textcolor{red}{$9\!\cdot\!10^{-5}\, x$}
            										&	$a=9\!\cdot\!10^{-5},\; b=10^{-7}$			&	\textcolor{red}{0.62}	\\
            			&				&	quadratic	&	0									&	$a=b=c=10^{-7}$								&	373.1	\\
            			&				&	exponential	&	$0.0001\; e^{0.0001\,x}$				&	$a=b=0.0001$								&	135.5	\\
            \midrule
            Cubic		&	Gurobi ILP
            			&	log			&	$25.39\,\ln x$						&	$a=25.39,\; b=10^{-6}$						&	230.7	\\
            			&				&	\textcolor{red}{\bf{linear}}
            										&	\textcolor{red}{$0.00791\, x$}
            										&	$a=0.00791,\; b=10^{-7}$						&	\textcolor{red}{27.58}	\\
            			&				&	quadratic	&	33.4									&	$a=b=10^{-7},\; c=33.4$						&	56.7	\\
            			&				&	exponential	&	$0.065\; e^{0.0001\,x}$				&	$a=0.0001,\; b=0.065$						&	230.5	\\
            \cmidrule{2-6}
            			&	Gurobi IQP
            			&	log			&	$226.5\,\ln x$						&	$a=226.5,\; b=10^{-6}$						&	2187.6	\\
            			&				&	linear		&	$0.127\, x$							&	$a=0.127,\; b=10^{-7}$						&	329.7	\\
            			&				&	\textcolor{red}{\bf{quadratic}}
            										&	\textcolor{red}{$7\cdot10^{-7}\,x^{2}+0.09\,x$}
            										&	$a=7\cdot10^{-7},\; b=0.09,\; c=10^{-7}$		&	\textcolor{red}{171.9}	\\
            			&				&	exponential	&	$19.7\; e^{0.0001\,x}$				&	$a=0.0001,\; b=19.7$						&	1723.2	\\
            \bottomrule
        \end{tabularx}
}%
    \caption{Curve fitting for the computation time with respect to system size for the different optimization problem types and solvers. The log function is $y=a \ln{x} + b$, the linear function is $y=a x + b$, the quadratic function is $y=a x^2 + b \cdot x + c$, and the exponential function is $y=e^{x a} \cdot b$. Here $x$ denotes the system size (specifically, the number of spins in the underlying Ising model) and the coefficients $a$, $b$, $c$ are fit using least squares optimization. The dependent variable which these curves are being fit to, for each system size $x$, is the raw compute time in the case of Gurobi computations and is the TTS in the case of simulated annealing (average compute time over the $10$ random instances for each problem size). The lowest RMSE function fit is highlighted as red text for each problem type and solver setting. }
    \label{table:curve_fitting}
\end{table}

All runtime results are reported on a large ensemble of Ising models. The Ising models are defined over a range of problem sizes starting at $100$ variables and going up to $160,000$ binary decision variables, with larger intervals as the problem size increases. $10$ random instances are generated for each problem size - the goal being to quantify whether there is a dramatic difference in average case, hardest case, and easiest case runtime as the number of variables increases. 

Figure \ref{fig:random_sampling_energies} shows two random sampling normalized energy distributions for these classes of Ising models both with and without cubic terms, as the system size is increased. The random energy samples were converted into approximation ratios by the formula 

\begin{equation}
\text{Approximation Ratio} = \frac{C_\text{Max} - C_e}{C_\text{Max} - C_\text{Min}}
\end{equation}

where $C_e$ is the energy of a single spin configuration, $C_\text{Max}$ is the global maximum energy, and $C_\text{Min}$ is the global minimum energy. The maximum and minimum energies are computed using Gurobi (for all remaining simulations the Gurobi computations are strictly minimization). 

Figure~\ref{fig:random_sampling_energies} is an initial and simple way to probe the energy spectrum of these problem instances - in particular we see that there is concentration of the random sampling energy near an approximation ratio of $0.5$ as the system size increases. This means that the problem instances are not trivially easy to solve in the sense that randomly sampling does not give a large approximation ratio as the system size grows.

Figure~\ref{fig:quadratic_Gurobi_time} shows the Gurobi runtime scaling as a function of system size for the problem instances with no higher-order terms, where we see not only very low absolute runtimes but also a scaling that appears to be approximately linear. For these Ising models with no cubic terms, since the binary IQP approach already works very well at finding the optimal solution, we do not additionally run the binary ILP method with Gurobi. 

Figure~\ref{fig:cubic_Gurobi_time} shows the Gurobi runtime scaling as a function of system size for the problem instances with cubic terms. Figure~\ref{fig:cubic_Gurobi_time} shows clearly that the order-reduction binary IQP method requires more computation time compared to the binary ILP method -- a potential reason for this could be the addition of auxiliary variables in the binary IQP case. Both Figure~\ref{fig:quadratic_Gurobi_time} and Figure~\ref{fig:cubic_Gurobi_time} show approximately linear time scaling. The average case runtime in Figures~\ref{fig:quadratic_Gurobi_time} and \ref{fig:cubic_Gurobi_time} does not substantially deviate from the largest or smallest runtimes, suggesting that these measured runtimes are reasonably consistent across different random polynomial coefficient choices.

\begin{figure}[ht!]
    \centering
    \includegraphics[width=0.49\linewidth]{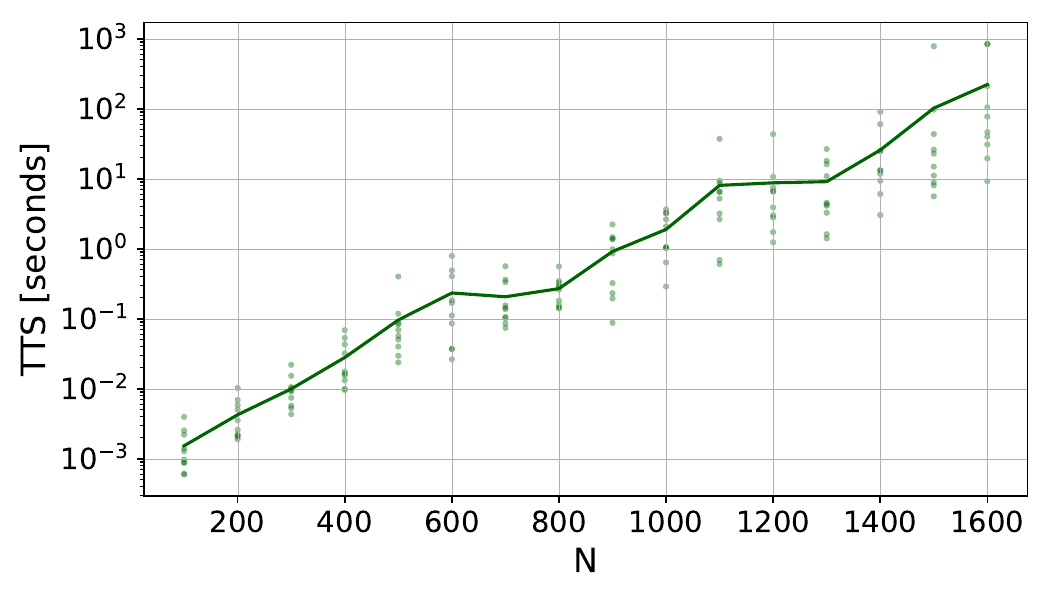}
    \includegraphics[width=0.49\linewidth]{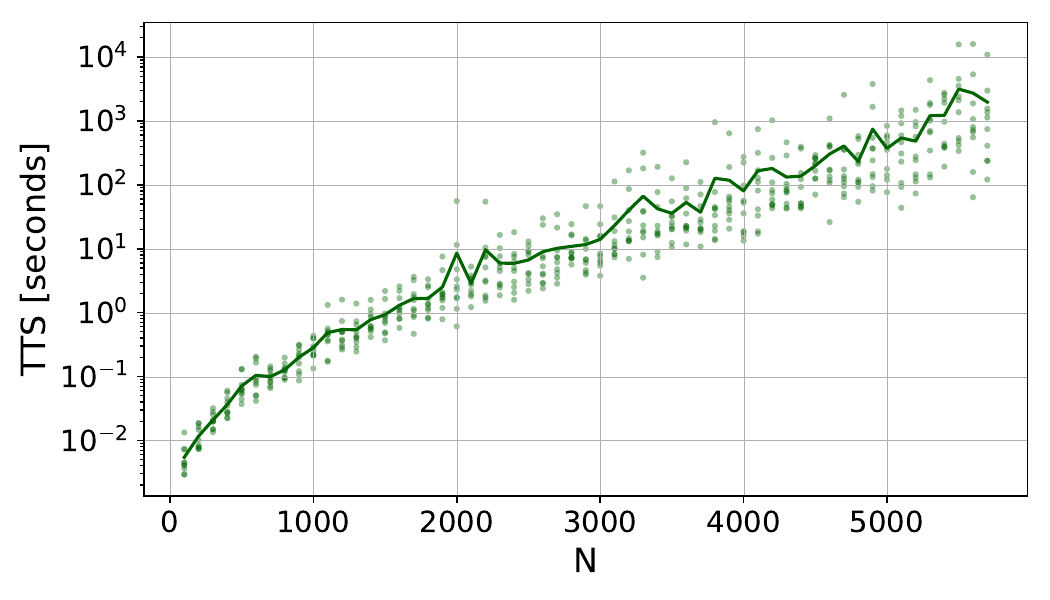}
    \caption{Simulated annealing TTS scaling (log scale y-axis) as a function of Ising model problem size (x-axis) for the instances containing no higher order terms. The x-axis (problem size) is exactly the number of nodes in the heavy-hex graph. $100$ Metropolis-Hastings updates (left) and $1000$ Metropolis-Hastings updates (right). At each system size N, $10$ different random-coefficient instances are solved and then the average time to solve those $10$ instances are shown as a function of $N$ by the connected green line. Here, $N$ refers specifically to the number of nodes in the underlying heavy-hex graph (in other words, not including any auxiliary variables produced from the order reduction process). TTS is measured using $50,000$ (independently generated) samples produced for each problem instance. }
    \label{fig:quadratic_SA_TTS}
\end{figure}

Figure~\ref{fig:quadratic_SA_TTS} shows time to solution scaling as a function of system size for the problem instances with no higher-order terms, with $100$ and $1000$ fixed variable sweep updates. Figure~\ref{fig:quadratic_SA_TTS} datapoints are only shown up until there is still finite sampling of an optimal solution, using the fixed number of samples computed for each instance. Figure~\ref{fig:quadratic_SA_TTS} appears to show exponential scaling. 

Table~\ref{table:curve_fitting} quantifies least squares fitting coefficients for four possible types of scaling functions namely logarithmic, linear, quadratic, and exponential. How well these fitted functions model the measured computation time data is quantified in terms of root-mean-square-error (RMSE). The lowest RMSE for all Gurobi runtime scaling are linear and quadratic function fits, for the problem instances solved using both problem encoding methods of ILP and IQP. The lowest RMSE for the Time to Solution at $99\%$ confidence (TTS$_{0.99}$) using simulated annealing sampling is the exponential function. The simulated annealing TTS exponential curve fits are weakly exponential, with exponents of $0.00898$ and $0.002$, for $100$ and $1000$ variable update sweeps respectively.

\begin{figure}[ht!]
    \centering
    \includegraphics[width=0.32\linewidth]{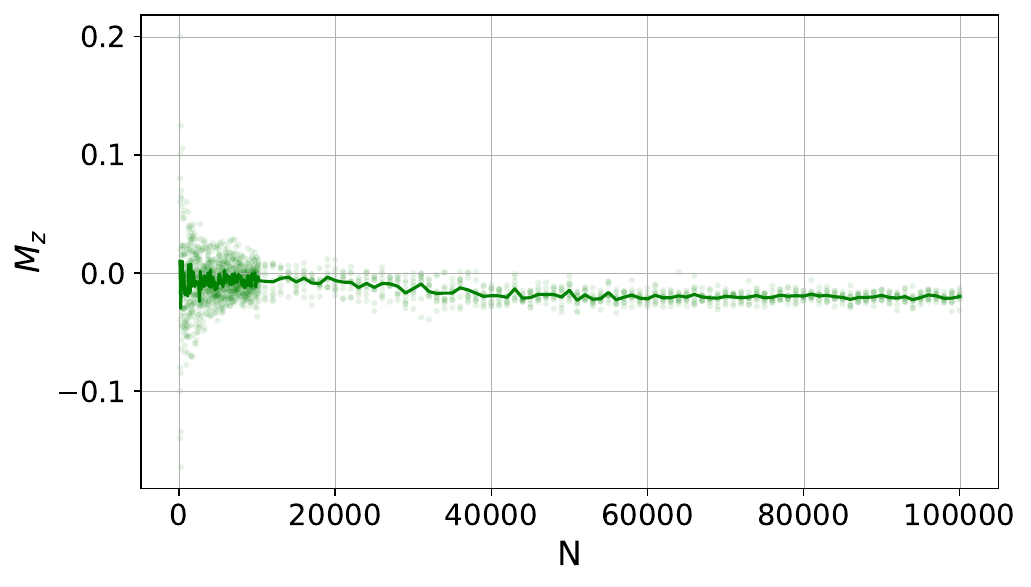}
    \includegraphics[width=0.32\linewidth]{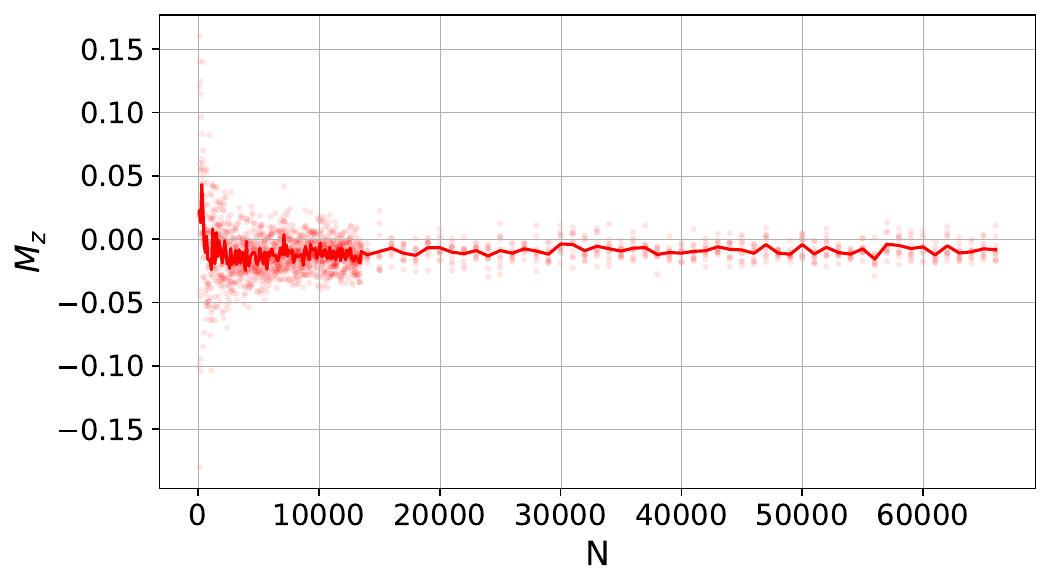}
    \includegraphics[width=0.32\linewidth]{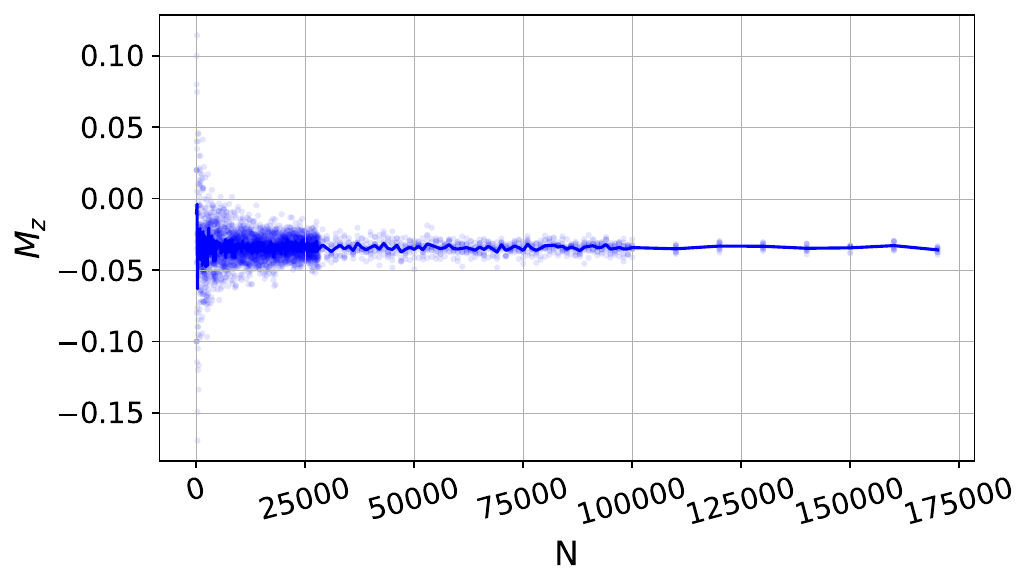}
    \caption{Optimal solution magnetization (y-axis) as a function of system size (x-axis), for the two cubic term problems (left and middle), and the quadratic term problems only (right). The average optimal solution magnetization is slightly below $0$ for all three sub-plots. 
    The magnetization of each (single) optimal solution found by Gurobi is shown, note however that this does not include exhaustive generate ground-state enumeration; this is magnetization of only a single degenerate ground-state for each problem instance (if the problem does have ground-state degeneracy). }
    \label{fig:magnetization}
\end{figure}

Having the capability to obtain verified optimal ground-states of these optimization problems, using Gurobi, means that we can examine some properties of these ground-states, such as net magnetization. Magnetization is a very simple magnetic order parameter, which serves as a simple observable test case that shows that we can use a commercial optimization solver to measure ground-state properties of this class of hardware-defined spin glass models. More complex ground-state order parameters could be extracted from these types of hardware-compatible Ising models, but we leave this to future work. Figure~\ref{fig:magnetization} shows the average $Z$-basis magnetization of the optimal solution found by Gurobi\footnote{Note that here we are measuring the net magnetization of the single optimal solution found, per problem instance, found by the Gurobi software.}, which shows a slight negative net magnetization, which becomes more concentrated as the system size increases, and is consistent for the problem instances with and without geometrically local cubic terms. Net magnetization of approximately zero is consistent with this class of Ising models, that being in general $\pm J$ models since we have approximately equal contributions from both antiferromagnetic and ferromagnetic interactions. In particular it means that the ground state does not have a very strongly preferred magnetic direction, but a slight ferromagnetic bias. Note that these plots do not show how consistent this property is for different degenerate ground-states of these systems. This net magnetization of the ground-state could be a property of this class of Ising models, due to their structure and coefficient distribution, or this could be a result of the algorithms used by the Gurobi software.

\section{Discussion and Conclusion}
\label{section:conclusion}

Primarily because this class of Ising models are incredibly sparse, and quite structured, they are relatively computationally easy for classical algorithms such as commercial mathematical programming software. There are ways that problems despite being very sparse can be computationally challenging, but at least for these classes of problems that this study examines this is not the case -- part of this could be due to the geometric locality of the higher order terms, and also potentially because of the low-coefficient precision of the problem instances. The fundamental motivation for these optimization problems is NISQ hardware-compatibility for a specific type of Quantum Processing Unit (QPU), which allows for scalable near-term quantum computer and algorithmic benchmarking of QAOA \cite{QAOA_QA_127, pelofske2023short, pelofske2023scalingwholechipqaoahigherorder}. However, in order to find optimization problems that are more computationally challenging for classical algorithms (and therefore a potential case where quantum algorithms could yield substantial improvements), promising optimization problem instances will need to include more complex interactions thus moving away from being purely hardware-compatible. Examples of more complex optimization types include longer range variable interactions and higher-weight high-order terms that are not necessarily geometrically local. The findings of the results in this study contextualizes where quantum algorithms performance for solving (whether approximately or exactly) combinatorial optimization problems need to get to in order to be competitive with state of the art classical algorithms, albeit for a specific type of NISQ processor hardware graph. More complexly connected optimization problems also mean that with current noisy quantum processors the degree of meaningful signal that can be extracted is not expected to be very high (as opposed to very hardware compatible optimization problems)~\cite{pelofske2023scalingwholechipqaoahigherorder, pelofske2023short, Preskill_2018, Andrist_2023, Harrigan_2021}. Therefore, in the near term, the central question is whether there exist interesting optimization problems that are both relatively NISQ hardware compatible, but also challenging for state of the art classical optimization solvers. It seems very likely that, if such problems exist, they will need to be on the scale of hundreds of decision variables since problems smaller than that are typically solved very efficiently by existing algorithms. 

Note that ground-state degeneracy is not quantified using any of the classical solvers used in this study, but this could be a topic for future study.

Importantly, the variable types that are considered in this study being discrete, and in particular spins, are well suited for digital quantum algorithms such as QAOA. And therefore, in future comparisons of quantum algorithms and classical algorithms, it will be necessary to compare against solvers that handle this variable type. However, it is worth noting that many natural types of optimization problems do not inherently have decision variable types that are binary or spins, and in such cases casting to something like an Ising model can be detrimental for the performance of linear programming solvers. In such cases, fair comparisons should be made against state of the art linear programming solvers by \emph{not} going through extensive problem re-formulation, or order reductions, when utilizing the mathematical optimization software.

\section*{Acknowledgments}
\label{sec:acknowledgments}
We thank Carleton Coffrin for discussions on integer programming and Gurobi. This work was supported by the U.S. Department of Energy through the Los Alamos National Laboratory. Los Alamos National Laboratory is operated by Triad National Security, LLC, for the National Nuclear Security Administration of U.S. Department of Energy (Contract No. 89233218CNA000001). Research presented in this article was supported by the NNSA's Advanced Simulation and Computing Beyond Moore's Law Program at Los Alamos National Laboratory. This research used resources provided by the Darwin testbed at Los Alamos National Laboratory (LANL) which is funded by the Computational Systems and Software Environments subprogram of LANL's Advanced Simulation and Computing program (NNSA/DOE). 

\noindent
LA-UR-24-32607

\bibliography{references}

\end{document}